\documentclass[11pt]{article}
\usepackage{amsmath}
\usepackage{amssymb}
\usepackage{eucal}
\usepackage{stmaryrd}
\usepackage[usenames]{color}

\begin{document}

\baselineskip 16pt

\title{On       $H_{\sigma}$-permutably
   embedded  and $H_{\sigma}$-subnormaly
   embedded subgroups  of finite groups  \thanks{Research is supported by
a NNSF grant of China (Grant \# 11371335) and Wu Wen-Tsun Key
 Laboratory of Mathematics of Chinese Academy of Sciences}   }
  
\author{Wenbin Guo,  Chi Zhang\\
{\small Department of Mathematics, University of Science and
Technology of China,}\\ {\small Hefei 230026, P. R. China}\\
{\small E-mail:
wbguo@ustc.edu.cn}\\ \\
{ Alexander  N. Skiba, Darya A. Sinitsa}\\
{\small Department of Mathematics,  Francisk Skorina Gomel State University,}\\
{\small Gomel 246019, Belarus}\\
{\small E-mail: alexander.skiba49@gmail.com}}

\date{}
\maketitle

\begin{abstract}
Let  $G$ be a finite group. Let  $\sigma =\{\sigma_{i} | i\in I\}$ be a
partition of the set of all primes $\Bbb{P}$ and $n$ an integer.
We write    $\sigma (n)
=\{\sigma_{i} |\sigma_{i}\cap \pi (n)\ne  \emptyset  \}$,
$\sigma (G) =\sigma (|G|)$.
A set  $ {\cal H}$ of subgroups of $G$ is said to be a  \emph{complete
Hall $\sigma
$-set} of $G$ if   every   member  of  ${\cal H}\setminus \{1\}$ is a Hall
 $\sigma _{i}$-subgroup of $G$ for some $\sigma _{i}$ and ${\cal H}$ contains
 exact one Hall  $\sigma
_{i}$-subgroup of $G$ for every  $\sigma _{i}\in  \sigma (G)$.
A  subgroup $A$ of  $G$ is called: (i) a \emph{$\sigma$-Hall
 subgroup} of $G$  if $\sigma (|A|) \cap  \sigma (|G:A|)=\emptyset$;
  (ii)     \emph{${\sigma}$-permutable}    in $G$  if $G$ possesses  a complete
Hall $\sigma$-set  ${\cal H}$  such that $AH^{x}=H^{x}A$  for all $H\in {\cal H}$
and all $x\in G$.
We say that a  subgroup $A$ of $G$ is  \emph{$H_{\sigma}$-permutably
 embedded}  in  $G$  if $A$ is a  ${\sigma}$-Hall subgroup of some
${\sigma}$-permutable    subgroup of $G$.

We study  finite groups $G$ having  an $H_{\sigma}$-permutably
 embedded  subgroup of order $|A|$ for each subgroup $A$ of $G$.
Some known results are generalized.

\end{abstract}


\let\thefootnote\thefootnoteorig

\section{Introduction}

Throughout this paper, all groups are finite and $G$ always denotes
a finite group.   Moreover,  $n$ is an integer,
 $\mathbb{P}$ is the set of all    primes, and if
  $\pi \subseteq  \Bbb{P}$, then   $\pi' =  \Bbb{P} \setminus \pi$.
  The symbol $\pi (n)$
denotes the
 set of all primes dividing $n$; as usual,  $\pi (G)=\pi (|G|)$, the set of all
 primes dividing the order of $G$.  We use $n_{\pi}$ to denote the  ${\pi}$-part of $n$,
 that is,  the largest ${\pi}$-number dividing   $n$; $n_{p}$ denotes the largest degree
 of $p$ dividing $n$.

In what follows, $\sigma =\{\sigma_{i} | i\in I\}$ is  some
partition of $\Bbb{P}$, that is,
$\Bbb{P}=\cup_{i\in I} \sigma_{i}$ and $\sigma_{i}\cap
\sigma_{j}= \emptyset  $ for all $i\ne j$;  $\Pi$ is  a
 subset of $\sigma$  and $\Pi'=\sigma\setminus
\Pi$.

Let    $\sigma (n)
=\{\sigma_{i} |\sigma_{i}\cap \pi (n)\ne  \emptyset  \}$ and
 $\sigma (G)
=\sigma (|G|)$. Then we say that  $G$ is 
\emph{$\sigma$-primary} \cite{1} if $G$ is a $\sigma _{i}$-group for some  $\sigma _{i}\in
 \sigma$.

A set  $ {\cal H}$ of subgroups of $G$ is said to be a  \emph{ complete
Hall $\sigma
$-set} of $G$  \cite{2, 3} if   every  member of  ${\cal H}\setminus \{1\}$ is a Hall
 $\sigma _{i}$-subgroup of $G$ for some $\sigma _{i}$ and ${\cal H}$ contains
 exact one Hall  $\sigma
_{i}$-subgroup of $G$ for every  $\sigma _{i}\in  \sigma (G)$.
We say that $G$ is  \emph{$\sigma
$-full} if $G$  possesses a complete
Hall $\sigma$-set.  Throughout this paper,  $G$ is always supposed  to be
a {$\sigma $-full  group.

A  subgroup $A$ of $G$ is called \cite{1}: (i) a \emph{$\sigma$-Hall
 subgroup} of $G$  if $\sigma (|A|) \cap  \sigma (|G:A|)=\emptyset$;
  (ii)
   \emph{${\sigma}$-subnormal}
  in $G$  if
there is a subgroup chain  $$A=A_{0} \leq A_{1} \leq \cdots \leq
A_{t}=G$$  such that  either $A_{i-1} \trianglelefteq A_{i}$ or
$A_{i}/(A_{i-1})_{A_{i}}$ is  ${\sigma}$-primary
  for all $i=1, \ldots ,t$;
  (iii)   \emph{${\sigma}$-quasinormal}  or    \emph{${\sigma}$-permutable}    in $G$  if $G$ possesses  a complete
Hall $\sigma$-set  ${\cal H}$  such that $AH^{x}=H^{x}A$  for all $H\in {\cal H}$
and all $x\in G$.  In particular, $A$ is called \emph{$S$-quasinormal}  or 
  \emph{$S$-permutable} in $G$ \cite{prod, GuoII} provided  
$AP=PA$ for all Sylow subgroups $P$ of $G$.

 {\bf Definition 1.1.}   We say that
 a subgroup $A$ of $G$ is     \emph{  $H_{\sigma}$-subnormally}
 (respectively \emph{$H_{\sigma}$-permutably},  $H_{\sigma}$-normally})
 embedded in  $G$   if $A$ is a  ${\sigma}$-Hall subgroup of some
 ${\sigma}$-subnormal   (respectively ${\sigma}$-permutable,  normal)  subgroup of $G$.

In the  special case, when $\sigma =\{\{2\}, \{3\}, \ldots  \}$,
 the definition of $H_{\sigma}$-normally embedded subgroups is
 equivalent to the concept of  Hall normally embedded subgroups in
\cite{2009},   the definition of $H_{\sigma}$-permutably  embedded subgroups is equivalent to
 the concept of  Hall $S$-quasinormally  embedded subgroups in \cite{clt}
and the definition  of $H_{\sigma}$-subnormally embedded subgroups is
 equivalent to the concept of  Hall subnormally embedded subgroups in
\cite{shirong}.

{\bf Example 1.2.}    (i) For any $\sigma$, all $\sigma$-Hall subgroups 
  and all ${\sigma}$-subnormal subgroups  of any group $S$
 are $H_{\sigma}$-subnormally embedded in $S$. Now, let  $G=(C_{7}\rtimes
C_{3}) \times A_{5}$, where  $C_{7}\rtimes
C_{3}$ is a non-abelian group of order 21 and  $A_{5}$ is the alternating
 group of degree 5, and  let $H=(C_{7}\rtimes
C_{3})  A$, where $A$ is a Sylow 2-subgroup of $A_{5}$.
 Let $\sigma  =\{\sigma _{1}, \sigma _{2}\}$,
 where
$\sigma _{1}=\{7\}$   and $\sigma _{2} =\{7\}'$.
Then $H$ is
 $\sigma$-subnormal in $G$ and $C_{3}A_{5}$ is a $\sigma$-Hall
 subgroup of $G$. In view of Lemma 2.1(1)(5) below,  the subgroup  $C_{3}A$ is neither
$\sigma$-subnormal in $G$ nor $H_{\sigma}$-normally embedded in $G$.

(ii)   For any $\sigma$, all $\sigma$-Hall subgroups and all $\sigma$-permutable   subgroups
 of any group $S$ are $H_{\sigma}$-permutably  embedded in  $S$.  Now,
 let
 $p > q > r$ be   primes, where  $r^{2}$ divides $q-1$. Let
 $\sigma =\{\sigma _{1}, \sigma _{2}  \}$,
 where
  $\sigma _{1}= \{q, r\}$ and
$\sigma _{2} =\{q, r\}'$.
 Let   $H=Q\rtimes R$
 be a  group of order $qr^{2}$, where  $C_{H}(Q)=Q$.  Let $P$ be  a simple
  ${\mathbb F}_{p}H$-module  which is     faithful  for $H$ 
 and  $G=P\rtimes H$. Let $R_{1}$ be a subgroup of $R$ of order $r$.
  Then the subgroup   $V=PR_{1}$ is $\sigma $-permutable  in $G$ and
 $R_{1}$ is a $\sigma$-Hall subgroup of $V$.  Hence $R_{1}$ is
$H_{\sigma}$-permutably  embedded in  $G$.
 It is also  clear that $G$ has
 no an  $S$-permutable   subgroup $W$ such that $R_{1}$ is a Hall subgroup of $W$, so
$R_{1}$ is neither  $H_{\sigma}$-normally  embedded nor $S$-permutably embedded  in  $G$.

  (iii) For any $\sigma$, all $\sigma$-Hall subgroups and all normal subgroups
 of any group $S$ are $H_{\sigma}$-normally  embedded in  $S$.   Now,  let
 $P$  be a simple  ${\mathbb F}_{11}(C_{7}\rtimes
C_{3})$-module
 which is     faithful  for $C_{7}\rtimes
C_{3}$. Let $G=(P\rtimes (C_{7}\rtimes C_{3})) \times A_{5}$.
 Let $\sigma  =\{\sigma _{1}, \sigma _{2}\}$,
 where
$\sigma _{1}=\{5, 7, 11\}$   and $\sigma _{2} =\{5, 7, 11\}'$.
Then the subgroup   $M= (P\rtimes C_{7}) \times A_{5}  $ is normal in $G$
  and a  subgroup $B$ of
$A_{5}$ of order 12  is  ${\sigma} $-Hall subgroup of $M$, so
$B$ is $H_{{\sigma} }$-normally  embedded in  $G$. Finally,  if
$\sigma  =\{\sigma _{1}, \sigma _{2}\}$,
 where
$\sigma _{1}=\{7\}$   and $\sigma _{2} =\{7\}'$, then
$B$ is not   $H_{{\sigma} }$-normally  embedded in  $G$.

Recall  that $G$ is  \emph{$\sigma$-nilpotent} \cite{4}
 if  $G=H_{1}\times \cdots \times H_{t}$ for some $\sigma$-primary groups
$ H_{1}, \ldots ,H_{t} $.
    The \emph{$\sigma$-nilpotent} residual $G^{{\frak{N}}_{\sigma}}$ of $G$
 is the intersection of all normal subgroups $N$ of $G$ with
 $\sigma$-nilpotent quotient $G/N$, $G^{\frak{N}}$ denotes  the nilpotent residual  of $G$.
It is clear that  every subgroup of a $\sigma$-nilpotent group $G$ is  
$\sigma$-permutable  and  $\sigma$-subnormal in $G$.

 {\bf Theorem  1.3.}
 {\sl Let ${\cal H}= \{1, H_{1}, \ldots H_{t} \}$ be a complete Hall
 $\sigma$-set of $G$ and   $D=G^{{\frak{N}}_{\sigma}}$.
 Then any two of the following conditions are equivalent:}

(i) {\sl   $G$  has   an $H_{\sigma}$-permutably
 embedded   subgroup of order $|A|$ for each subgroup $A$ of $G$. }

 (ii) {\sl $D$ is a complemented cyclic of square-free order  subgroup of 
$G$ 
 and $|\sigma _{i} \cap \pi (G)|=1$ for all $i$ such that
 $\sigma _{i} \cap \pi (D)\ne \emptyset$.}

(iii) {\sl For each  set $\{A_{1}, \ldots ,A_{t}\}$, where $A_{i}$ is a   subgroup
  (respectively  normal subgroup)
 of $H_{i}$ for all $i=1, \ldots ,t$,
 $G$  has an $H_{\sigma}$-permutably  embedded
 (respectively $H_{\sigma}$-normally  embedded) subgroup
  of order  $|A_{1}| \cdots  |A_{t}|$.}

Let $\mathfrak{F}$ be a class of groups. A subgroup $H$ of $G$ is
 said to be an $\mathfrak{F}$-covering subgroup of $G$ \cite[VI, Definition
 7.8]{hupp} if $H\in \mathfrak{F}$ and  for every subgroup $E$ of $G$ such
 that $H\leq E$ and $E/N\in \mathfrak{F}$ it follows that $E=NH$.
  We say that a  subgroup  $H$ of $G$ is a
  \emph{$\sigma$-Carter subgroup} of $G$ if $H$ is an
  $\mathfrak{N}_{\sigma}$-covering subgroup of $G$, where
 $\mathfrak{N}_{\sigma}$ is the class of all $\sigma$-nilpotent groups.

A group $G$ is said to have a \emph{Sylow tower} if $G$  has a
 normal series $1=G_{0} < G_{1} < \cdots < G_{t-1} < G_{t}=G$,  where 
 $|G_{i}/G_{i-1}|$ is the order of
 some Sylow subgroup of $G$ for each  $i\in \{1, \ldots , t\}$.  
A chief factor of $G$  is said to be 
 \emph{$\sigma$-central} (in $G$) \cite{1} if the semidirect product
 $(H/K) \rtimes (G/C_{G}(H/K))$  is $\sigma$-primary.
  Otherwise,   $H/K$ is called  \emph{$\sigma$-eccentric} (in $G$).

We  say that $G$ is a \emph{$H\sigma E$-group} if the following conditions are hold:
(i)   $G=D\rtimes M$, where
  $D=G^{{\frak{N}}_{\sigma}}$ is a  $\sigma$-Hall subgroup of $G$ 
  and $|\sigma (D)|=|\pi (D)|$; (ii) 
   $D$ has a Sylow tower and every chief factor of
 $G$ below $D$ is $\sigma$-eccentric; (iii)  $M$
 acts irreducibly on every $M$-invariant Sylow subgroup of $D$.

We do not still know the structure of a  group $G$ having  a  $H_{\sigma}$-subnormally
 embedded   subgroup of order $|A|$ for each subgroup $A$ of $G$.
Nevertheless, the following fact is true.

{\bf Theorem  1.4.} {\sl  Any two of the following conditions are
equivalent:}

(i) {\sl Every subgroup of $G$ is  $H_{\sigma}$-subnormally
 embedded   in $G$.}

(ii) {\sl  Every ${\sigma}$-subnormal subgroup $H$ of $G$
 is an  $H\sigma E$-group of the form $H=D\rtimes M$, where  
$D=H^{{\frak{N}}_{\sigma}}$   and 
 $M$ is a   $\sigma $-Carter subgroup of $H$.  }

 (iii) {\sl   Every ${\sigma}$-subnormal subgroup  of $G$
 is an  $H\sigma E$-group.}

Now,  let us consider some corollaries of Theorems 1.3 and 1.4.  First note that
 since a nilpotent group $G$ possesses a normal subgroup of order $n$ for each integer
 $n$ dividing $|G|$,  in the   case when $\sigma =\{\{2\}, \{3\}, \ldots
\}$, Theorem 1.3
 covers Theorem 11 in \cite{2009}, Theorem 2.7 in \cite{shirong} and Theorems
 3.1 and 3.2 in \cite{clt}.

From Theorem 1.3 we  also get the following  result.

 {\bf Corollary  1.5.}
 {\sl Suppose that $G$ possesses a  complete Hall
 $\sigma$-set ${\cal H}= \{1, H_{1}, \ldots ,H_{t} \}$
 such that $H_{i}$ is nilpotent for all $i=1, \ldots ,t$.
   Then   $G$  has   an $H_{\sigma}$-normally
 embedded   subgroup of order $|H|$ for each subgroup $H$ of $G$ if and only if the nilpotent
residual $D=G^{\frak{N}}$ of $G$ is cyclic of square-free order and
 $|\sigma _{i} \cap \pi (G)|=1$ for all $i$ such that
 $\sigma _{i} \cap \pi (D)\ne \emptyset$.}

In the   case when $\sigma =\{\{2\}, \{3\}, \ldots  \}$ we get from
Corollary 1.5
 the following known result.

 {\bf Corollary   1.6} (Ballester-Bolinches,  Qiao \cite{cltII}).
 {\sl $G$  has   a  Hall  normally
 embedded   subgroup of order $|H|$ for each subgroup $H$ of $G$ if and only if
 the nilpotent
residual $G^{\frak{N}}$ of $G$ is cyclic of square-free order.}

On the basis of Theorems 1.3 and 1.4 we prove also  the next two theorems.

{\bf Theorem   1.7.} {\sl  Any two of the following conditions are
equivalent:}

(i) {\sl Every subgroup of $G$ is  $H_{\sigma}$-normally
 embedded   in $G$.}

(ii) {\sl    $G=D
 \rtimes M$
 is a  $H\sigma E$-group, where  $D$
 is a cyclic group  of  square-free order and    $M$ is a     Dedekind  group.  }

 (iii) {\sl  $G=D\rtimes M$, where
  $D$ is a  $\sigma$-Hall cyclic  subgroup of $G$
 of  square-free order
 with $|\sigma (D)|=|\pi (D)|$ and    $M$ is a     Dedekind  group.}

In the   case when $\sigma =\{\{2\}, \{3\}, \ldots  \}$ we get from
Theorem  1.7
 the following  known  result.

{\bf Corollary 1.8} (Li, Liu  \cite{shirong}).  {\sl Every subgroup of
 $G$ is a Hall  normally embedded  subgroup of   $G$ if and only if
 $G=D\rtimes M$, where
  $D=G^{\frak{N}}$ is a cyclic Hall subgroup of $G$ of
 square-free order
and $M$ is a  Degekind  group.}

{\bf Theorem  1.9.} {\sl  Any two of the following conditions are
equivalent:}

(i) {\sl Every subgroup of $G$ is  $H_{\sigma}$-permutably
 embedded   in $G$.}

(ii) {\sl    $G=D \rtimes M$
 is a  $H\sigma E$-group, where 
$D=G^{{\frak{N}}_{\sigma}}$
 is a  cyclic group  of  square-free order.  }

 (iii) {\sl  $G=D\rtimes M$, where
  $D$ is a  $\sigma$-Hall cyclic  subgroup of $G$
 of  square-free order
 with $|\sigma (D)|=|\pi (D)|$ and   $M$  is ${\sigma}$-nilpotent}

{\bf Corollary 1.10.} {\sl Every subgroup of $G$
 is   a Hall    $S$-quasinormally  embedded  subgroup  of $G$ if and only if
 $G=D\rtimes M$, where
  $D=G^{{\frak{N}}}$ is a cyclic Hall subgroup of $G$ of
 square-free order
and $M$ is a   Carter subgroup of $G$.}

In conclusion of this section, consider the following example.

{\bf Example 1.11.}   Let $5 < p_{1} < p_{2}  < \cdots < p_{n}$ be
a set of primes and $p$ a prime such that either $p > p_{n}$ or $p$ divides $p_{i}-1 $ for all
$i=1, \ldots , n$. Let
 $A$ be a
 group of order $p$
and $P_{i}$ a simple  ${\mathbb F}_{p_{i}}A$-module
 which is     faithful  for $A$. Let $L_{i}=P_{i}\rtimes A$
 and $$B=(\ldots ((L_{1}\Yup L_{2})\Yup L_{3})  \Yup  \cdots )\Yup L_{n}$$
 (see \cite[p. 50]{hupp}). We can assume without loss of generality that $L_{i}\leq B$
 for all
$i=1, \ldots , n$. Let
  $G=B \times A_{5}$, where $A_{5}$ is the
 alternating group of degree 5,   and  let 
$\sigma$ be  a partition of $\mathbb{P}$   such that for some
different indices $i, j, i_{1},\ldots   , i_{n}\in I$  we have $p\in \sigma _{i}$,
 $\{2, 3, 5\}\subseteq
 \sigma _{j}$ and  $p_{k}\in  \sigma _{i_{k}}$ for all $k=1, \ldots , n$.
Then  $$D=P_{1}P_{2} \cdots  P_{n}=G^{{\frak{N}}_{\sigma}}$$ is a
 $\sigma$-Hall subgroup of $G$ and $G=D\rtimes (A\times A_{5})$.

We show that  every subnormal subgroup $H$ of $G$  satisfies Condition (ii) in Theorem 1.4.
If  $H^{{\frak{N}}_{\sigma}}= 1$, it is evident. Hence we can assume 
without loss of generality  $A\leq H$ since
 every $p'$-subgroup of $G$ is $\sigma$-nilpotent. 
 But then $$H=(H\cap D)\rtimes (A\times (H\cap A_{5}))$$ by Lemma 2.1(4), where $H\cap D$ is
a  normal $\sigma$-Hall subgroup of $H$  and $M=A\times
 (H\cap A_{5})$ is a $\sigma$-nilpotent
subgroup of $H$. Moreover,  $H\cap A_{5}$  induces on every non-identity Sylow subgroup of
 $H\cap D$ a non-trivial irreducible
group of automorphisms. Therefore  $H^{{\frak{N}}_{\sigma}}=H\cap D$
 and $|\sigma (H^{{\frak{N}}_{\sigma}})|=|\pi
(H^{{\frak{N}}_{\sigma}})|$. It is  also clear that $M$ is a
 $\sigma$-Carter subgroup of $H$ and  every chief factor of $H$ below
 $H^{{\frak{N}}_{\sigma}}$ is
   $\sigma$-eccentric in $H$.   Thus $G$
satisfies Condition (ii) in Theorem 1.4, and so  every subgroup $H$   of $G$ is
 $H_{\sigma}$-subnormally  embedded  in  $G$. On the other hand, the
subgroup $DAC_{2}$, where $C_{2}$ is a subgroup of order 2 of $G$,
is not Hall subnormally  embedded  in  $G$ since $C_{2}$
 is not a Sylow subgroup of any subnormal subgroup of $G$.

Finally, if  $p$ divides $p_{i}-1 $ for all
$i=1, \ldots , n$, then $|P_{i}|=p_{i}$  for all
$i=1, \ldots , n$,  so $G$ satisfies Condition (ii) in Theorem 1.9 and
hence  satisfies Condition (ii) in Theorem 1.3.

\section{Basic lemmas}

An integer $n$ is called a \emph{$\Pi$-number} if $\sigma (n) \subseteq
\Pi$.   A  subgroup $H$ of $G$ is called a \emph{Hall $\Pi$-subgroup} of
$G$ \cite{1} if $|H|$ is a  $\Pi$-number and $|G:H|$ is a $\Pi'$-number.

{\bf Lemma 2.1} (See Lemma 2.6 in \cite{1}). {\sl Let  $A$,  $K$ and
 $N$ be subgroups of  $G$, where   $A$
is $\sigma$-subnormal in $G$ and $N$ is normal in $G$.  }

(1) {\sl $A\cap K$    is  $\sigma$-subnormal in
$K$}.

(2) {\sl If $K$ is  $\sigma$-subnormal in $G$, then $A\cap K$
  and  $\langle A, K\rangle$  are     $\sigma$-subnormal in
$G$.}

(3) {\sl $AN/N$ is
$\sigma$-subnormal in $G/N$. }

(4) {\sl If $H\ne 1 $ is a Hall $\Pi$-subgroup of $G$  and $A$ is not  a
 $\Pi'$-group, then $A\cap H\ne 1$ is
 a Hall $\Pi$-subgroup of $A$. }

(5) {\sl If  $|G:A|$ is a $\sigma_{i}$-number, 
 then  $O^{{\sigma_{i}}}(A)= O^{{\sigma _{i}}}(G)$.}

(6) {\sl  If $V/N$ is  a $\sigma$-subnormal subgroup of  $G/N$, then $V$ is
$\sigma$-subnormal in $G$. }

(7) {\sl If $K$
is a $\sigma$-subnormal subgroup of  $A$, then  $K$
is  $\sigma$-subnormal in   $G$}.

 A group  $G$ is said to be \emph{$\sigma$-soluble }\cite{1}  if every chief factor of $G$
 is $\sigma$-primary.

{\bf Lemma 2.2} (See Lemmas 2.8 and  3.2 and Theorems B and C in \cite{1}).  {\sl 
 Let  $A$,  $K$ and
 $N$ be subgroups of  $G$, where  $A$
is $\sigma$-permutable in $G$ and $N$ is normal in $G$.}

(1)  {\sl   $AN/N$ is $\sigma$-permutable   in $G/N$.}

(2)  {\sl If  $G$ is    $\sigma$-soluble, then
 $A\cap K$ is $\sigma$-permutable in $K$.}

(3)  {\sl If  $N\leq K$, $K/N$
 is  $\sigma$-permutable in $G/N$ and $G$ is   $\sigma$-soluble,
 then $K$ is  $\sigma$-permutable in $G$.}

(4)   {\sl $A$
is $\sigma$-subnormal   in $G$.}

(5)   {\sl If  $G$ is   $\sigma$-soluble  and  $K$
is $\sigma$-permutable   in $G$,  then  $K\cap A$
is $\sigma$-permutable   in $G$.}

 {\bf Lemma  2.3.} {\sl Let  $H$ be a normal subgroup of $G$. If
 $H/H\cap \Phi (G)$ is a $\Pi$-group, then $H$ has a a Hall $\Pi$-subgroup,
 say  $E$, and $E$ is normal in $G$.  Hence, if $H/H\cap \Phi (G)$ is $\sigma$-nilpotent,
 then $H$ is $\sigma$-nilpotent.}

{\bf Proof.} Let $D=O_{\Pi'}(H)$. Then, since $H\cap \Phi (G)$ is
nilpotent,  $D$ is a Hall  $\Pi'$-subgroup of $H$. Hence by  the Schur-Zassenhaus
theorem, $H$ has a Hall $\Pi$-subgroup, say  $E$.  It is clear
that $H$ is $\pi'$-soluble where $\pi' =\cup _{\sigma _{i}\in \Pi'}\sigma _{i}$,
 so any two Hall $\Pi$-subgroups of $H$ are conjugate.
By the Frattini argument,  $$G=HN_{G}(E)=(E(H\cap \Phi (G)))N_{G}(E)=N_{G}(E).$$   Therefore
 $E$   is
normal in $G$.  The lemma is proved.

{\bf Lemma 2.4.} {\sl If  every chief factor of $G$ below $D=G^{{\frak{N}}_{\sigma}}$
 is cyclic, then $D$ is nilpotent.}

{\bf Proof.} Assume that this
 is false and let $G$ be a counterexample of minimal order.
Let $R$ be a minimal normal subgroup of $G$. Then from the $G$-isomorphism
 $$D/D\cap R\simeq DR/R=  (G/R)^{{\frak{N}}_{\sigma}}$$ we know that  every
 chief factor of $G/R$ below $DR/R$ is cyclic, so the choice of $G$ implies that
 $D/D\cap R\simeq DR/R$ is nilpotent. Hence $R\leq D$ and $R$ is the unique minimal
 normal subgroup of $G$. In view of Lemma 2.3,
 $R\nleq \Phi (G)$ and  so $R=C_{R}(R)$ by  \cite[A, 15.2]{DH}. But by hypothesis,
 $|R|$ is a prime, hence
  $G/R=G/C_{G}(R)$ is cyclic, so $G$ is supersoluble and so
 $G^{{\frak{N}}_{\sigma}}$ is nilpotent     since $G^{{\frak{N}}_{\sigma}}\leq
 G^{{\frak{N}}}$. The lemma is proved.

The  following lemma is evident.

{\bf Lemma 2.5.}  {\sl   The class of all $\sigma$-soluble groups
 is closed under taking direct
products, homomorphic images and  subgroups. Moreover, any extension of the
 $\sigma$-soluble group by a $\sigma$-soluble group is
 a $\sigma$-soluble group as well.    }

Let $A$, $B$ and $R$ be subgroups of $G$. Then $A$ is said to \emph{$R$-permute}
with  $B$ \cite{GS8} if for some $x\in R$ we have $AB^{x}=B^{x}A$.

 If $G$ has a complete Hall $\sigma$-set ${\cal H}=\{1, H_{1}, \ldots ,
H_{t} \}$   such that $H_{i}H_{j}=H_{j}H_{i}$ for
all $i, j,$ then we say  that $\{H_{1}, \ldots ,
H_{t} \}$    is a \emph{$\sigma$-basis}
of $G$.

{\bf Lemma 2.6} (See Theorems A and B in \cite{2}).  {\sl Assume that 
 $G$ is $\sigma$-soluble.  }

(i) {\sl  $G$ has  a $\sigma$-basis
$\{H_{1}, \ldots , H_{t} \}$ such that
for each $i\ne j$ every Sylow subgroup of $H_{i}$  $G$-permutes with every
Sylow subgroup of $H_{j}$. }

(ii) {\sl For any $\Pi$, the following hold:
 $G$ has a  Hall $\Pi$-subgroup $E$,  every  $\Pi$-subgroup of $G$ is
contained in some conjugate of $E$ and $E$   $G$-permutes with every
Sylow subgroup of $G$.   }

{\bf Lemma 2.7.} {\sl Let  $H$,  $E$ and $R$ be subgroups of  $G$.
 Suppose that   $H$
is $H_\sigma$-subnormally  embedded  in $G$ and $R$ is normal in $G$.}

(1) {\sl If    $H\leq E$, then $H$ is    $H_{\sigma}$-subnormally embedded
 in $E$.}

(2) {\sl  $HR/R$ is    $H_{\sigma}$-subnormally embedded
 in   $G/R$.}

(3) {\sl If $S$  is  a  ${\sigma}$-subnormal subgroup of $G$,  then $H\cap S$ is
    $H_{\sigma}$-subnormally embedded  in $G$.}

(4)  {\sl If $|G:H|$ is ${\sigma}$-primary, then $H$ is either a ${\sigma}$-Hall
subgroup of $G$ or
    ${\sigma}$-subnormal  in $G$. }

{\bf Proof.}  Let $V$ be a ${\sigma}$-subnormal subgroup of $G$
 such that $H$ is a $\sigma$-Hall subgroup of $V$.

 (1) This assertion is a corollary of  Lemma 2.1(1).

(2)  In view of Lemma 2.1(3),    $VR/R$ is
      ${\sigma}$-subnormal   subgroup of $G/R$. It is also clear  that $HR/R$ is
 a $\sigma$-Hall subgroup of $VR/R$.
Hence  $HR/R$ is $H_{\sigma}$-subnormally embedded   in  $G/R$.

(3)  By Lemma 2.1(1)(2),   $V\cap S$ is   ${\sigma}$-subnormal both  in
$V  $ and  in  $G$ and so  $H\cap (V\cap S)=H\cap S$ is
 a  ${\sigma}$-Hall subgroup of   $V\cap S$ by Lemma 2.1(4),  as required.

(4)  Assume that $H$ is not  ${\sigma}$-subnormal in $G$. Then  $H < V$. By hypothesis,
  $|G:H|$ is
${\sigma}$-primary, say $|G:H|$ is a ${\sigma _{i}}$-number. Then $|V:H|$
is a ${\sigma _{i}}$-number.  But $H$ is a ${\sigma}$-Hall subgroup  of    $V$.
Hence   $H$ is a ${\sigma}$-Hall subgroup  of $G$.

The lemma is proved.

{\bf Lemma  2.8.} {\sl Let $H$   be a    ${\sigma}$-subnormal subgroup of
 a    ${\sigma}$-soluble group $G$.  If   $|G:H|$ is a $\sigma _{i}$-number and
$B$ is a $\sigma _{i}$-complement of $H$, then $G=HN_{G}(B)$.  }

{\bf Proof.}  Assume that this lemma is false and let $G$ be a
counterexample of minimal order. Then $H < G$, so $G$ has a proper subgroup $M$ such that $H \leq M$,
$|G:M_{G}|$ is a $\sigma _{i}$-number and $H$ is  ${\sigma}$-subnormal  in $M$. The choice
 of $G$ implies that $M=H N_{M}(B)$. On the other hand, clearly that $B$ is
 a  $\sigma_{i}$-complement of $M_{G}$.  Since $G$ is  ${\sigma}$-soluble, Lemma 2.6 and
the Frattini argument imply that $$G=M_{G}N_{G}(B)=MN_{G}(B)=H N_{M}(B)N_{G}(B)=HN_{G}(B).$$
 The lemma is proved.

The following lemma is well-known (see for example \cite[Lemma 3.29]{Shem-Sk} or \cite[1.10.10]{Guo}).

{\bf Lemma  2.9.} {\sl Let  $H/K$ be an abelian  chief factor of $G$
   and  $V$  a maximal
subgroup of $G$ such that $K\leq V$ and $HV=G$. Then $$G/V_{G}\simeq
(H/K)\rtimes (G/C_{G}(H/K)).$$}

\section{Proofs of the results}

 {\bf Proof of Theorem 1.3.} Without   loss of generality we may assume
that $H_{i}$ is a $\sigma _{i}$-group for all $i=1, \ldots , t$.

(i), (iii)  $\Rightarrow$ (ii)  Assume that this   is false.
 Then $D\ne 1$ and so $t > 1$.

First we prove the following claim. 

(*) {\sl If   $p\in \sigma _{i}\cap \pi (G)$, then  $G$ has  a
   $\sigma$-permutable  subgroup $E$   with $|E|=|G|_{\sigma _{i}'}p$. }

We can assume without loss of generality that $i=1$. In fact, to prove 
Claim (*), we consistently build  the $\sigma$-permutable subgroups   $E_{2}, \ldots , E_{t}$
 such that   $|H_{2}|\cdots
|H_{j}|$  divides    $|E_{j}|$ and $|E_{j}|_{\sigma_{1}}=p$ for all $j=2,
\ldots ,t$.

By hypothesis,  $G$ has an  $H_{\sigma}$-permutably  embedded
  subgroup $X$ of order $p$. Let $V$ be a  $\sigma$-permutable  subgroup
 of $G$ such  that  $X$ is a  $\sigma$-Hall
subgroup of   $V$. Then  $|V|_{\sigma _{1}}=p$ and $G$ has  a complete Hall $\sigma$-set
 $\{1, K_{1},
 \ldots ,K_{t} \}$, where $K_{i}$ is a $\sigma _{i}$-group for all $i=1, \ldots , t$,
  such that $VK_{i}=K_{i}V$ for all $i=1, \ldots , t$.  Let $W=VK_{2}$. 
Then   $|W|_{\sigma _{1}}=p$. 

Next we show that  there is an  $H_{\sigma}$-permutably  embedded subgroup $Y$ of $G$
  such that $|Y|=|W|$.  It is enough   
to consider the case when  Condition (iii) holds.   
  Let $A_{1}$ be a subgroup of 
$H_{1}$ of 
order  $p$, $A_{2}= H_{2}$ and $A_{i}= H_{i}\cap V$ for all $i > 2$. 
Then $$|A_{2}|= |H_{2}|=|K_{2}|.$$    
 On the other hand, $V\cap K_{i}$ and  $V\cap H_{i}$ are Hall $\sigma _{i}$-subgroups 
of $V$ by Lemmas 2.1(4)  and 2.2(4) and so $|V\cap K_{i}|= |V\cap H_{i}|$.
Also, 
      for every    $ i > 2$ we have 
$$|W:V\cap K_{i}|=|VK_{2}: V\cap K_{i}|=|V||K_{2}|:|V\cap K_{2}|| 
V\cap K_{i}|$$  is a $\sigma _{i}'$-number and hence  $V\cap K_{i}=W\cap K_{i}$ is a 
Hall $\sigma _{i}$-subgroup  of $W$. Therefore,  $$|W|= p|H_{2}||V\cap 
H_{3}|\ldots   |V\cap 
H_{t}|$$ and so  $G$  has an $H_{\sigma}$-permutably  embedded
  subgroup $Y$ 
  of order  $$|W|=|A_{1}| \cdots  |A_{t}|$$ by hypothesis.

Let $E_{2}$ be a  $\sigma$-permutable  subgroup
 of $G$ such  that  $Y$ is a  $\sigma$-Hall
subgroup of   $E_{2}$.   Then $|H_{2}|=|K_{2}|$   divides
  $|E_{2}|$ and $|E_{2}|_{\sigma _{1}}=p$.    Now, arguing by induction, 
 assume that  $G$
 has
a $\sigma$-permutable  subgroup  $E_{t-1}$ such that   $|H_{2}|\cdots
|H_{t-1}|$  divides    $|E_{t-1}|$ and $|E_{t-1}|_{\sigma _{1}}=p$.
 Then  for some Hall  $\sigma _{t}$-group  $L$ we have  
$E_{t-1}L=LE_{t-1}$, and  if $E_{t}=E_{t-1}L$, then $|E_{t}|=|G|_{\sigma _{1}'}p$ and
 $E_{t}$ clearly is $\sigma$-permutable in $G$, as required.

Now, let   $p\in \sigma_{i}\cap \pi (D)$ and let  $P$ be a Sylow $p$-subgroup of $D$. 
Then, by Claim (*),    $G$ possesses  a
$\sigma$-permutable  subgroup  $E$ such that   $|E|=|G|_{\sigma _{i}'}p$.
Lemma 2.2(4) implies that $E$ is
$\sigma$-subnormal in $G$.  Let $j\neq i$. Then $H_{j}^{x}\cap E$ 
is a Hall $\sigma _{j}$-subgroup of $E$ by
Lemma 2.1(4), so $|E:H_{j}^{x}\cap E|$ is a  ${\sigma _{j}}'$-number.
 But $|H_{j}^{x}|$ divides $|E|$ and hence 
  $|H_{j}^{x}|$ divides $|H_{j}\cap E|$.
  Thererore  $H_{j}^{x}\leq E$  for all $x\in G$. Thus $H_{j}^{G}\leq E$ 
and so  
$G/E_{G}$ is a   $\sigma _{i}$-group. Hence  $D\leq E_{G}\leq E$, so   $|P|= p$.
Therefore  $G$ is supersoluble by \cite[IV,  2.9]{hupp} and so every chief 
factor of $G$ below $D$
 is cyclic. Hence  $D$ is nilpotent by Lemma 2.4, so $D$ is  cyclic of 
square-free order. Hence $D$ is complemented in $G$ by Theorem 11.8 in  \cite{shem}. 

Finally, assume that $|\sigma _{i} \cap \pi 
(G)| > 1$ and let  
 $q\in \sigma _{i} \cap \pi (G)\setminus \{p\}$.  Then   $G$ possesses  a
$\sigma$-permutable  subgroup  $F$ such that   $|F|=|G|_{\sigma _{i}'}q$. 
Then 
$D\leq F_{G}\leq F. $ Therefore $D\leq E\cap F$ and so $p$ does not divide $|D|$.  
    This contradiction completes
the proof of the implications (i) $\Rightarrow$ (ii) and (iii) $\Rightarrow$ (ii).

(ii) $\Rightarrow$ (iii)
 First we show that
 for every $i$ and for every
     subgroup  (respectively normal subgroup) $A_{i}$ of $H_{i}$,
   there is an $H_{\sigma}$-permutably embedded
 (respectively $H_{\sigma}$-normally embedded)
     subgroup $E_{i}$ of $G$ such
that $|E_{i}|=|A_{i}||G|_{\sigma _{i}'}$.   
 Since  $G$  evidently  is $\sigma$-soluble, it
 has a  $\sigma_{i}$-complement  $E$ by Lemma 2.6.    Therefore, it is enough to
 consider the case when $A_{i}\ne 1$ since every $\sigma$-Hall subgroup of $G$
  is an   $H_{\sigma}$-normally
     embedded in  $G$.

First suppose that  $D\leq E$. Then  $E/D$ is normal in $G$
 since $G/D$ is $\sigma$-nilpotent. Therefore $$(E/D)\times (A_{i}D/D)=EA_{i}/D$$
is $\sigma$-permutable   (respectively normal)  in $$G/D=  (E/D) \times (H_{i}D/D).$$
 Hence
$E_{i}=EA_{i}$  is  $\sigma$-permutable   (respectively normal)  in $G$ by Lemma
2.2(3)  and
 $|E_{i}|=|A_{i}||G|_{\sigma _{i}'}$.

Now suppose that $D\nleq E$. Then $D\cap H_{i}\ne 1$, so $H_{i}$ is a $p$-group
 for some prime $p$ since for each $\sigma _{i} \in \sigma (D)$ we have
  $|\sigma _{i}\cap \pi (G)|=1$ by hypothesis. Hence $H_{2} $ has    a normal subgroup
  $A$ such that $D_{p}\leq A$ and $|A|=|A_{i}|$, where $D_{p}$ is a Sylow $p$-subgroup of $D$.
 Then  $D\leq  AE$. Moreover,
 $$AE/D=(DA/D)\times (ED/D)$$
 since $ED/D$ is a Hall $\sigma _{i}'$-subgroup of $G/D$.
 Therefore $E_{i}=AE$ is ${\sigma}$-permutable    (respectively normal)
  in $G$ by Lemma 2.2(3) and  $|E_{i}|=|A_{i}||G|_{\sigma _{i}'}$.

 Let  $E=E_{1}\cap \cdots \cap E_{t}$.
 Then $|E|=|A_{1}|\ldots |A_{t}|$ since
 $$(|G:E_{i}|, |G:E_{j}|)=1$$  for all $i\ne j$. Note that $E_{i}$ is either a
 $\sigma$-Hall subgroup of $G$ or  $\sigma $-permutable (respectively normal)  in $G$.
 Indeed, let $V$ be a $\sigma $-permutable   (respectively normal)  subgroup of
 $G$ such that $E_{i}$ is a
 $\sigma$-Hall subgroup of $V$.
 Assume that $E_{i}$ is not  $\sigma $-permutable (respectively not normal)   in $G$.
 Then  $E_{i}
 < V$.
 Since   $|G:E_{i}|$ is
${\sigma _{i}}$-number,  $|V:E_{i}|$
is a ${\sigma _{i}}$-number.  But $E_{i}$ is a ${\sigma}$-Hall subgroup  of    $V$.
Hence   $E_{i}=V$ is a ${\sigma}$-Hall subgroup  of $G$.

  Assume that $E_{1}, \ldots ,
 E_{r}$ are   $\sigma $-permutable (respectively normal) in $G$ and $E_{i}$
 is a $\sigma$-Hall subgroup of $G$ for
 all $i > r$. Then $E^{0}=E_{1} \cap \cdots \cap E_{r}$ is  $\sigma
$-permutable  (respectively normal)  in $G$ by Lemma 2.2(5) and
 $E^{*}=E_{r+1} \cap \cdots \cap E_{t}$ is a $\sigma$-Hall subgroup of $G$. Now,
  $E=E^{0}\cap E^{*}$ is  a $\sigma$-Hall subgroup of $E^{0}$ by Lemmas 2.1(4) and 2.2(4),
 so $E$ is
  $H_{\sigma}$-permutably  (respectively $H_{\sigma}$-normally) embedded in
 $G$.
 Hence  (ii) $\Rightarrow$ (iii).

(ii) $\Rightarrow$ (i) Since  $G$ is $\sigma$-soluble, $H$ is
$\sigma$-soluble. Hence $H$ has a  $\sigma$-basis $\{L_{1}, \ldots , L_{r}\}$  such that
$L_{i} \leq H_{i}$
 for
all $i=1, \ldots, r$ by Lemma 2.6. Therefore from the implication  (ii)
$\Rightarrow$ (iii) we get that $G$ has an $H_{\sigma}$-permutably  embedded
  subgroup
  of order  $|L_{1}| \cdots  |L_{r}|=|H|$.

The theorem is proved.

{\bf Proof of Theorem 1.4.}   (i) $\Rightarrow$ (ii)    Assume that this
 is false and let $G$ be a
counterexample of minimal order. Then some  $\sigma$-subnormal subgroup
$V$ of $G$ is not an $H\sigma E$-group.
 Moreover,  $D=G^{{\frak{N}}_{\sigma}} \ne 1$, so $|\sigma (G)| > 1$.

(1) {\sl   Condition (ii) is true   on
every proper section $H/K$ of $G$, that is, $K\ne 1$ or $H\ne G$. Hence  $V=G$}
 (This directly follows from Lemma 2.7(1)(2) and the choice of $G$).

(2) {\sl $G$ is ${\sigma}$-soluble.}

In view of Claim (1) and Lemma 2.5, it is enough to show that $G$ is not simple. Assume
that this is false. Then 1 is the only proper $\sigma$-subnormal subgroup of $G$  since
$|\sigma (G)| > 1$.
Hence every  subgroup of $G$ is a $\sigma$-Hall subgroup of $G$.  
Therefore  for a Sylow $p$-subgroup $P$ of $G$, where $p$  is the smallest prime divisor of
$|G|$, we have $|P|=p$ and so  $|G|=p$ by \cite[IV, 2.8]{hupp}.
This contradiction shows that we have (2).

(3) {\sl If 
$|G:H|$ is a   $\sigma _{i}$-number  and $H$ is not a $\sigma$-Hall
subgroup of $G$, then  $H$ is  $\sigma$-subnormal in $G$ and a
$\sigma _{i}$-complement $E$ of $H$ is normal in $G$} (This follows from Lemmas 2.7(4) 
  and 2.8). 
  
(4) {\sl  $D$ is a Hall subgroup of $G$. Hence $D$ has a complement $M$ in $G$.}

 Suppose
that this is false and let $P$ be a  Sylow $p$-subgroup of $D$ such
that $1 < P < G_{p}$, where $G_{p}\in \text{Syl}_{p}(G)$.  We can assume 
without loss of generality that $G_{p}\leq H_{1}$. Let $R$ be a minimal normal 
subgroup of $G$ contained in $D$.

 Since
 $D$ is  soluble by Claim (2),   $R$ is a $q$-group    for some prime   
$q$. Moreover, 
$D/R=(G/R)^{\mathfrak{N}_{\sigma}}$  is a Hall subgroup of $G/R$ by
Claim (1) and Proposition 2.2.8 in \cite{15}.  Suppose that  $PR/R \ne 1$. Then  $PR/R \in \text{Syl}_{p}(G/R)$. 
If $q\ne p$, then    $P \in \text{Syl}_{p}(G)$. This contradicts the fact 
that $P < G_{p}$.  Hence $q=p$, so $R\leq P$ and therefiore $P/R \in \text{Syl}_{p}(G/R)$.
 It follows that 
$P \in \text{Syl}_{p}(G)$. This contradiction shows that  $PR/R=1$, which implies that 
  $R=P$ is a Sylow $p$-subgroup of 
$D$. Therefore $R$ is a  unique minimal normal 
subgroup of $G$ contained in $D$. It is also clear  that a $p$-complement of $D$ is a Hall 
subgroup of $G$.

 Now we  show that  $R\nleq \Phi (G)$. Indeed, assume that  $R\leq \Phi (G)$.
  Then  $D\ne R$ by Lemma 2.3 since $D=G^{{\frak{N}}_{\sigma}}$.
On the other hand,    $D/R$ is a  $p'$-group. Hence $O_{p'}(D)\ne 1$  
 by Lemma 2.3. But   $O_{p'}(D)$ is characteristic in $D$ and so 
 it is 
normal $G$. Therefore $G$ has a minimal normal subgroup $L$ such that  $L\ne R$ and $L\leq D$.
This contradiction shows that  $R\nleq \Phi (G)$.

 Let $S$ be a maximal subgroup of $G$ such that $RS=G$.  Then $|G:S|$ is a
 $p$-number. Hence, since $R$ is not a Sylow $p$-subgroup of $G$, $p$ divides $|S|$. Then  $S$ is not  a Hall
subgroup of $G$ and so $S$ is not  a $\sigma$-Hall
subgroup of $G$. Therefore    $S$ is    $\sigma$-subnormal in $G$ by Claim (3) and  so
 $G/S_{G}$ is a  $\sigma _{i}$-group, which implies that $$R\leq D\leq S_{G}\leq S$$  and so
 $G=RS=S$.
This contradiction completes the proof of (4).

(5)  {\sl If  $M\leq   E < G$, then  $E$ 
 is not  $\sigma$-permutable in $G$ and  so  $E$ a 
 $\sigma$-Hall subgroup of $G$.}

Assume that $E$ is    $\sigma$-permutable in $G$.  Then $E$ is    $\sigma$-subnormal in $G$ by Lemma 2.2(4). 
Then there is a subgroup chain  $$E=E_{0} \leq
E_{1} \leq \cdots \leq E_{r}=G$$ such that
either $E_{i-1}$   is normal in $E_{i}$ 
  or $E_{i}/(E_{i-1})_{E_{i}}$ is  $\sigma $-primary  for all $i=1, \ldots , r$. 
  Let   $V=E_{r-1}$.
  We can assume without loss of generality that $V\ne G$.  Therefore, 
since $G$ is $\sigma$-soluble by Claim (2), for some $\sigma$-primary 
chief factor $G/W$ of $G$ we have $E\leq V\leq W$.  Also we have  $D\leq   W$ and  so 
$G=DE\leq W$, a contradiction.   Hence $E$ is not    
$\sigma$-permutable in $G$.

  By hypothesis, $G$ has a $\sigma$-permutable subgroup $S$ such that
$E$ is a $\sigma$-Hall subgroup of $S$.  But then $S=G$, by the above argument, so
  $E$ is a $\sigma$-Hall subgroup of $G$. In particular, $M$ is a 
  $\sigma$-Hall subgroup of $G$ and so $D$ is a $\sigma$-Hall subgroup of 
$G$.

(6) {\sl $D$ is soluble,  $|\sigma (D)|=|\pi (D)|$ and    $M$
 acts irreducibly on every $M$-invariant Sylow subgroup of $D$}.

Let $p\in \sigma _{i}\in \sigma (D)$. Lemma 2.6
   and Claims (2) and  (4) imply that for some Sylow $p$-subgroup $P$ of $G$ we have
 $PM=MP$. Moreover,    $MP$ is a $\sigma$-Hall subgroup of $G$ by Claim (5). 
 Hence  $|\sigma _{i} \cap \pi (G)|=1$ for all $i$ such that
 $\sigma _{i} \cap \pi (D)\ne \emptyset$ and  so $|\sigma (D)|=|\pi (D)|$.  
 Therefore, since $D$ is soluble by Claim (2),
 $M$     acts irreducibly on every $M$-invariant Sylow subgroup of $D$
by  Claim (5).

(7)  {\sl $M$ is  a $\sigma$-Carter subgroup of $G$. }

Let $R$ be a minimal normal subgroup of $G$ contained in
 $D$ and $E$ a subgroup of $G$ containing $M$. We need to show
 that $E=E^{{\frak{N}}_{\sigma}}M$. Claim (1) implies that $RM/R$  is
  a $\sigma$-Carter subgroup of $G/R$, so $$ER/R=(ER/R)^{{\frak{N}}_{\sigma}}(RM/R).$$
Hence $ER=E^{{\frak{N}}_{\sigma}}MR$ since
 $$(ER/R)^{{\frak{N}}_{\sigma}}=E^{{\frak{N}}_{\sigma}}R/R.$$
 Claim (6) implies that $R$ is a $p$-group for some prime $p$.
 Claims (4), (5) and (6) imply that $R$, $E$ and $E^{{\frak{N}}_{\sigma}}M$ are
 $\sigma$-Hall subgroups of $G$. Therefore, if $R\nleq E$, then $E$ and
 $E^{{\frak{N}}_{\sigma}}M$ are    Hall $p'$-subgroups of
 $ER=E^{{\frak{N}}_{\sigma}}MR$, so  $E=E^{{\frak{N}}_{\sigma}}M$.
Finally, assume that $R\leq E$ but $R  \nleq E^{{\frak{N}}_{\sigma}}M$.
 Then $R\cap E^{{\frak{N}}_{\sigma}}=1$. On the other hand,
 since $$DE/D\simeq E/D\cap E$$ is $\sigma$-nilpotent, $E^{{\frak{N}}_{\sigma}} \leq D$
 and so $M\cap E^{{\frak{N}}_{\sigma}}=1$.
Therefore $$E^{{\frak{N}}_{\sigma}}\cap RM=
(E^{{\frak{N}}_{\sigma}}\cap R)(E^{{\frak{N}}_{\sigma}}\cap M)=1.$$ Then
 $$E/E^{{\frak{N}}_{\sigma}}=  E^{{\frak{N}}_{\sigma}}MR/E^{{\frak{N}}_{\sigma}}\simeq MR$$
 is $\sigma$-nilpotent. Hence $M\leq C_{G}(R)$. Suppose that  $
C_{G}(R)   < G$ and let $ C_{G}(R)\leq W <G$, where $G/W$ is a chief
factor of $G$. Claim (2) implies that $G/W$ is $\sigma$-primary, so $D\leq
W$. But then $G=DM\leq W < G$, a contradiction. Therefore $ C_{G}(R)=G$,
that is,  $R\leq Z(G)$.   Let $V$ be a complement to $R$ in $D$.
 Then $V$ is a Hall normal
 subgroup of $D$, so it is characteristic in $D$. Hence $V$ is normal in $G$ and
 $G/V\simeq RM$ is $\sigma$-nilpotent, so $D\leq V < D$. This contradiction completes
 the proof of (7).

(8) {\sl $D$ possesses a Sylow tower.}

Let $R$ be a minimal normal subgroup of $G$ contained in $D$. Then $R$ is
a $p$-group
 for some prime $p$ by Claim (6).  Moreover, the Frattini argument implies that for 
some Sylow $p$-subgroup $P$ of $D$  we have $M\leq N_{G}(P)$ and so $R=P$ 
since $M$ acts irreducible on $P$ by Claim (6). On the other hand, by 
Claim (1),  $D/R$   possesses
 a Sylow tower. Hence we have (8).

(9) {\sl Every chief factor of
 $G$ below $D$ is $\sigma$-eccentric.}

Let $H/K$ be  a chief factor of
 $G$ below $D$. Then $H/K$ is a $p$-group for some prime $p$ since $D$ is soluble by Claim (6).
By the Fratiini argument, there exist a Sylow
$p$-subgroup $P$ and   a $p$-complement
 $E$ of $D$  such that $M\leq N_{G}(P)$   and $M\leq N_{G}(E)$.
 Then $M\leq N_G(P\cap K)$  and $M\leq N_G(P\cap H)$. Hence $P\cap K=1$ 
and $P\cap H=P$ by Claim (6), so $H=K\rtimes P$.     
 Let $V=EM$. Then $K\leq V$ and $HV=G$, so  $V$ is a maximal
subgroup of $G$. Hence   $$G/V_{G}\simeq
(H/K)\rtimes (G/C_{G}(H/K))$$ by Lemma 2.9. Therefore, if $H/K$ is
$\sigma$-central in $G$, then $D\leq V_{G}$, which is impossible since
 evidently $p$ does not divide $|V|$.  Thus we have (9).

From Claims (4)--(9) it follows that $G$
 is a $H\sigma E$-group, contrary to our assumption on $G=V$.   Hence
(i) $\Rightarrow$ (ii).

(ii) $\Rightarrow$ (iii) This implication is evident.

(iii) $\Rightarrow$ (i)    By hypothesis,  $G=D\rtimes M$, where
$D=G^{{\frak{N}}_{\sigma}}$ is a  $\sigma$-Hall subgroup of $G$,  $|\sigma (D)|=|\pi (D)|$ and 
 $M$
 acts irreducibly on every $M$-invariant Sylow subgroup of $D$.

(*) {\sl Every 
subgroup $A$ of $G$ containing $M$ is a $\sigma$-Hall subgroup of $G$. }

    Let $D_{0}=D\cap A$. Then $A=D_{0}\rtimes M$ and 
  $D_{0}\ne 1$. Let $p\in \pi (D_{0})$. The Frattini argument and Lemma 2.6 imply  that 
for some Sylow    $p$-subgroup  $P_{0}$  of $D_{0}$ and  some Sylow $p$-subgroup $P$ of $D$ we 
have $M\leq  N_{G}(P_{0})$,  $M\leq  N_{G}(P)$ and $P_{0}M\leq PM$. 
Hence, since $M$
 acts irreducibly on every $M$-invariant Sylow subgroup of $D$,  
$P_{0}=P$.   Therefore every Sylow subgroup of $A$ is a Sylow subgroup of 
$G$. Hence  $A$ is a $\sigma$-Hall subgroup of $G$  since  $|\sigma 
(D)|=|\pi (D)|$  and $M$ is a $\sigma$-Hall subgroup of $G$.

Now, let $A$ be a  subgroup of $G$. First assume that $DA < G$. By Lemma 2.1(6), $DA$ is  $\sigma$-subnormal
in $G$. Therefore every $\sigma$-subnormal subgroup  of $DA$ is also
$\sigma$-subnormal in $G$. Hence  Condition (iii)  holds for $DA$, so  $A$
 is $H_{\sigma}$-subnormally embedded  in  $DA$ by induction. But then  $A$ is
  $H_{\sigma}$-subnormally embedded  in  $G$ by Lemma 2.1(7).

Finally, suppose that  $DA=G$.  Then, since   $G$  is $\sigma$-soluble, for
 some $x$ we have  $M\leq A^{x}$ by Lemma 2.6. Hence  $A^{x}$ is a   ${\sigma}$-Hall
 subgroup of  $G$ by Claim (*), so $A^{x}$ is an $H_\sigma$-subnormally embedded
subgroup of $G$. But then $A$ is an $H_\sigma$-subnormally embedded
subgroup of $G$.   Therefore   the implication (iii) $\Rightarrow$ (i) is
true.

  The theorem is proved.

{\bf Proof  of Theorem 1.9.}  (i) $\Rightarrow$ (ii)
    This follows from Lemma 2.2(4) and Theorems 1.3 and 1.4.

(ii) $\Rightarrow$ (iii)  This implication is evident.

(iii) $\Rightarrow$ (i)     Let $A$ be any subgroup of $G$. Then  $DA$
 is $\sigma$-permutable in $G$ by Lemma 2.2(3) since $G$ is
 $\sigma$-soluble.
 On the other hand, since $|\sigma (D)|=|\pi (D)|$ and $D$ is a cyclic $\sigma$-Hall
 subgroup of $G$ of square-free order,
$A$ is a $\sigma$-Hall subgroup of $DA$.
 Hence $A$ is  $H_\sigma$-permutably embedded in $G$.
  Therefore   the implication (iii) $\Rightarrow$ (i) is
true.

The theorem is proved.

{\bf Proof of Theorem  1.7. }
(i) $\Rightarrow$ (ii)  In view of Theorem 1.9, it is enough to show that if
 $D\leq L\leq G$ and $L$ is a   $\sigma$-Hall subgroup of some normal
 subgroup $V$ of $G$, then $L$ is normal in $G$. But since $G/D$ is
 $\sigma$-nilpotent,
  $L/D$ is $\sigma$-subnormal in $G/D$, so $L$ is  $\sigma$-subnormal in $G$ by Lemma
 2.1(6).
 Hence  $L$ is  $\sigma$-subnormal in $V$ by Lemma 2.1(1). But then
   $ L$ is a  normal in $V$ by Lemma 2.1(4) and so $L$ is a characteristic subgroup
 of $V$. It follows that $L$ is normal in $G$.

(ii) $\Rightarrow$ (iii)  This implication is evident.

 (iii) $\Rightarrow$ (i)  See the proof of the implication (iii) $\Rightarrow$ (i)  in
 Theorem 1.9.

The theorem is proved.

\end{document}